[1994/12/01]
\documentclass[draft]{article}
\usepackage{amsfonts}

\title{\bf Pseudo-hyperbolic distance and $n$-best rational approximation in $\mathbb{H}^2$ space
\thanks{This work is supported by NSFC 11501426}}

\author{Yan-Bo Wang\thanks{Faculty
of Mathematics and Computer Science, WuHan Textile University, WuHan, China,
E-mail:{\sl wybwtu@163.com}}, \ Tao Qian\thanks{Corresponding Author: Macao Center for Mathematical Sciences, Macao University of Science and Technology, Macao, China SAR, E-mail:{\sl tqian@must.edu.mo}}}
\usepackage{amsmath,amsthm}

\chardef\bslash=`\\ 
\newcommand{\ntt}{\normalfont\ttfamily}

\theoremstyle{definition}

\theoremstyle{remark}

\newcommand{\eval}[2][\right]{\relax
  \ifx#1\right\relax \left.\fi#2#1\rvert}

\begin{document}
\maketitle \markboth{Sample paper for the
{\protect\ntt\lowercase{amsmath}} package} {Sample paper for the
{\protect\ntt\lowercase{amsmath}} package}
\renewcommand{\sectionmark}[1]{}
{\bf Abstract.} Through reducing the problem to rational orthogonal system (Takenaka-Malmquist system), this note gives a proof for existence of $n$-best rational approximation to functions in the Hardy $\mathbb{H}^{2}(\mathbb{D})$ space by using pseudohyperbolic distance. \\
 \vspace{3mm} {\bf Key Words:} Schwarz Lemma, T-M System, $n$-Best Rational Approximation
\bigskip

{\bf AMS Subject Classifications:}  42A50, 32A30, 32A35, 46J15

\section{Introduction}
This paper concerns rational approximations in the Hardy space
\[ \mathbb{H}^{2}(\mathbb{D})=\{f:\mathbb{D}\to \mathbb{C} \ :\ f=\sum_{k=1}^\infty c_kz^k\ {\rm with} \sum_{k=1}^{\infty}|c_k|^2<\infty\},\]
where $\mathbb{D}$ is the open unit disc.
Among other equivalent definitions of the norm of the Hardy space, we adopt the one in terms of the non-tangential boundary limits of the functions
\[ \langle f,g\rangle =\frac{1}{2\pi}\int_0^{2\pi}f(e^{it})\overline{g}(e^{it})dt.\]
For a given positive integer $n,$ an ordered pair of polynomials $(p,q)$ is called an $n$-\emph{admissible pair} if $p$ and $q$ are co-prime, $q\ne 0,$ the zeros of $q$ are in the exterior of $\mathbb{D},$ and the degrees of $p$ and $q$ are both less than $n+1.$ The $n$-best rational approximation problem in the Hardy-$\mathbb{H}^{2}(\mathbb{D})$ space (abbreviated as $\mathbb{H}^2$) is stated as follows: For $f\in \mathbb{H}^{2}(\mathbb{D}),$ find an admissible pair $(\tilde{p},\tilde{q})$ such that
\begin{eqnarray}\label{problem}
\|f-\frac{\tilde{p}}{\tilde{q}}\|_{\mathbb{H}^2}=\min \{\|f-\frac{p}{q}\|_{\mathbb{H}^2} \ :\ (p,q) \ {\rm is\ an} \ n-{\rm admissible\ pair}\}.
\end{eqnarray}
For a Hardy space function $f$, denoting its non-tangential boundary limit also by $f$ if without confusion, there holds $\|f\|_{\mathbb{H}^2}=\|f\|_{L^2(\partial\mathbb{D})}.$ In below we will abbreviate both $\|\cdot\|_{\mathbb{H}^2}$ and $\|\cdot\|_{L^2(\partial\mathbb{D})}$ as $\|\cdot \|.$ This problem and closely related ones have been long studied. The existence of a solution to (\ref{problem}) under Chebyshev norm with prescribed poles was first proved in \cite{J.L.W1}, and existence, uniqueness of best rational approximation under $\|.\|_{p}$ norms were given in \cite{J.L.W2}. Due to the great interest, researchers have given alternative proofs based on different methods, of which some were related to algorithms to find a solution\cite{LB1}\cite{LB2}. The basic methodology before our approach introduced in\cite{Q1}\cite{MQ}, was to parameterize the problem by the coefficients of the denominator $q,$ and the related optimal numerator is obtained through orthogonalization \cite{J.L.W2}\cite{BN1}. Our approach to the problem is via Szeg\"o kernel approximation to functions in the space. By this approach we directly find the best suitable poles and so to find the best denominator.
The Szeg\"o kernel approach is as described now. Denote by $k_{a}$ the Szeg\"o kernel of the Hardy $\mathbb{H}^{2}(\mathbb{D})$ space, where
 \[ k_a(z)=\frac{1}{1-\bar{a}z},\]
 and by
\begin{eqnarray}
\mathfrak{D}:=\{e_{a}(z)=\frac{\sqrt{1-|a|^{2}}}{1-\bar{a}z}, a\in\mathbb{D}\},
\end{eqnarray}
  the collection of normalized Szeg\"o kernels. The function $k_a$ is the reproducing kernel in $\mathbb{H}^{2}(\mathbb{D})$.
If a sequence $\{a_{n}\}_{n=1}^{+\infty}\subset\mathbb{D}$ is given, where multiplicity is allowed, through the Gram-Schmidt orthogonalization process (G-S) on $\{e_{a_{n}}\}_{n=1}^{+\infty}$,
we can get an orthonormal
 system, called the rational orthogonal system or Takenaka-Malmquist system (T-M system),
\begin{eqnarray}
\mathfrak{E}:=\{E_{a_1,a_2,...,a_n}(z)=\frac{\sqrt{1-|a_{n}|^{2}}}{1-\bar{a}_{n}z}\prod_{k=1}^{n-1}\frac{z-a_{k}}{1-\bar{a}_{k}z}, n=1,2...\}
\end{eqnarray}
where $E_n=E_{a_1,a_2,...,a_n}=e_{a_n}B_{a_1,a_2,..,a_{n-1}},$ being the product of a normalized Szeg\"o kernel and an $(n-1)$-order Blaschke product. For any $n$ complex numbers $c_1,\cdots,c_n,$ the form
\[ \sum_{k=1}^n c_kE_k\]
is called an $n$-Blaschke form. If $c_n$ is none zero, then it is called an $n$-non-degenerated Blaschke form. For $f\in \mathbb{H}^2,$ there is an associated $n$-Blaschke form
\[ \sum_{k=1}^n \langle f,E_k\rangle E_k.\]
From the Hilbert space theory
\[ \|f-\sum_{k=1}^n \langle f,E_k\rangle E_k\|=\min\limits_{(c_1,\cdots,c_n)} \{\|f-\sum_{k=1}^n c_kE_k\|\}
\]
and
\[ \|f-\sum_{k=1}^n \langle f,E_k\rangle E_k\|=0\]
if and only if $f\in {\rm Span}\{E_k\}_{k=1}^n.$

There is correspondingly an $n$-best Blaschke form approximation problem: Find a set of $n$ parameters $a_1,\cdots,a_n,$ all in $\mathbb{D},$ such that
\begin{eqnarray}\label{nbest}& &\|f-\sum_{k=1}^n\langle f,E_{\!a_1,\cdots,a_k\!}\rangle E_{a_1,\cdots,a_k}\|\!\nonumber \\
&=&\!\inf \{\|f-\sum_{k=1}^n\langle f,E_{b_1,\cdots,b_k}\rangle E_{\!b_1,\cdots,b_k\!}\|: \!b_1,\cdots,b_n\!\ {\rm are\ all\ in }\ \mathbb{D}\}.\end{eqnarray}

Below we will refer (\ref{nbest}) as $\lq\lq$ $n$-best Blaschke form approximation ". For the connection between the $n$-best rational and the $n$-best Blaschke form functions there is the following observation: If there exists $a_k=0$ among an $n+1$ sequence $a_1,\cdots,a_{n+1},$ where the multiplicity is allowed, then  non-degenerate $(n+1)$-Blaschke form made from $a_1,\cdots,a_{n+1}$ are of the form $p/q,$ where $(p,q)$ is an $n$-admissible pair \cite{Q2}.

Based on the above observation, the $n$-best rational approximation problem and the $n$-best Blaschke approximation problem are essentially the same. In fact, to get an $n$-best rational approximation to $f$ one can, instead, solve the $n$-Blaschke form problem for $f-c_0,$ where $c_0$ is the $0$-th Fourier coefficient of $f.$ In the wide notion of sparse representation by linear combinations of the dictionary words in a Hilbert space with a dictionary, the problem in terms of $n$-Blaschke form seems to be more essential and natural, as well as more general. Our main result is as follows.\\

\noindent{\bf Theorem 1}
For any $f\in \mathbb{H}^2$ and any positive integer $n$, if $f$ is not identical with
an $m$-Blaschke form, $m<n,$ then there exists a non-degenerate $n$-best approximation to $f.$\\

A mathematical algorithm for finding a solution to the $n$-best Blaschke approximation has yet been an open problem. There, however, exist several proofs for existence of a solution that were mostly associated with the goal of finding an algorithm. In this note we provide a new proof for the existence by surprisingly using pseudohyperbolic distance. This methodology may lead a new way to treat similar problems in general reproducing kernel Hilbert functional spaces in which the Hardy space methods are unadaptable.

The pseudohyperbolic distance on $\mathbb{D}$ is defined by
\begin{eqnarray}
 \rho(z_0,z) &=&\left |\frac{z-z_0}{1-\bar{z}_{0}z}\right |.
\end{eqnarray}
If analytic function $g(z)$ is defined from $\mathbb{D}$ to $\overline{\mathbb{D}}$ , Schwarz lemma shows that\cite{J.B.G}
\begin{eqnarray}\label{schwarz lemma}
  \rho(g(z_0),g(z)) &\leq& \rho(z_0,z), \quad z_0\neq z
\end{eqnarray}
and
\begin{eqnarray}
  |g^{\prime}(z)|(1-|z|^2) &\leq& 1-|g(z)|^2.
\end{eqnarray}
In the sequel, we denote $f_{b_1,b_2,...,b_n}(z)=f(z)-\sum\limits_{k=1}^{n}\langle f,E_{b_1,b_2,...,b_k}\rangle E_{b_1,b_2,...,b_k}(z).$
 We note that
\[ \|f-\sum_{k=1}^n \langle f,E_{b_1,b_2,...,b_k}\rangle E_{b_1,b_2,...,b_k}\|^2=
\|f\|^2-\sum_{k=1}^n|\langle f,E_{b_1,b_2,...,b_k}\rangle|^2.\]
Hence, to attain
\[ \inf \|f-\sum_{k=1}^n \langle f,E_{b_1,b_2,...,b_k}\rangle E_{b_1,b_2,...,b_k}\|^2\]
is equivalent with to attain
\begin{eqnarray}\label{sup} \sup \sum_{k=1}^n|\langle f,E_{b_1,b_2,...,b_k}\rangle|^2.\end{eqnarray}
We will denote the orthogonal projection of $f$ into the linear subspace $X$ by
$P_{X}f.$ The projection into the subspace as orthogonal complement of $X$ is denoted $Q_Xf=(I-P_X)f.$ In the case $X={\rm Span}\{e_{a_1},\cdots,e_{a_n}\}$, they will be simply denoted as
$P_{{a_1},\cdots,{a_n}}$ and $Q_{{a_1},\cdots,{a_n}}.$ It is recognized that $Q_{{a_1},\cdots,{a_k}}$ is the Gram-Schmidt (G-S) process operator, and
 \[ E_{a_1,a_2,...,a_k}=Q_{a_1,a_2,...,a_{k-1}}(e_{a_k}).\]
 We note in this notation $f_{a_1,a_2,...,a_k}=Q_{{a_1},\cdots,{a_k}}f,$ and, owing to the self-adjoint property of projection operators and the orthogonality gained from the G-S process, for $a_l$ among $a_1,\cdots,a_k,$
\begin{eqnarray*}
 f_{a_1,a_2,...,a_k}(a_l)&=&\langle Q_{{a_1},\cdots,{a_k}}f, k_{a_l}\rangle\\
 &=&\langle f, Q_{{a_1},\cdots,{a_k}}k_{a_l}\rangle\\
 &=&\langle f, (I-P_{a_1,a_2,...,a_k})k_{a_l}\rangle\\
 &=& 0.\end{eqnarray*}

\section{Proof of Theorem}

We will first prove the following \\

\noindent{\bf Lemma 1.}
 Let $f$ be a Hardy space function with analytic continuation to $\overline{\mathbb{D}}$ with $\|f\|_{\mathbb{H}^\infty}\leq M$, then for any $k$-tuple  $a_1,\cdots,a_k\in \mathbb{D},$ there holds
$\|f_{a_1,\cdots,a_k}\|_{\mathbb{H}^\infty}\leq 3^kM.$\\

\noindent{\bf Proof.} When $k=1,$ by using the reproducing kernel property of $k_{a_1}(z)=\frac{1}{(1-\overline{a}_1z)},$ we have
\begin{eqnarray*}
|f_{a_1}(z)|&=&|f(z)-\langle f,E_{a_1}\rangle E_{a_1}(z)|\\
&\leq& |f(z)|+|f(a_1)|(1-|a_1|^2)\frac{1}{1-|a_1|}\\
&\leq& 3M.\end{eqnarray*}
Now treat the general $k>1$ case. For arbitrary $z\in \mathbb{D},$ due to the orthogonality and the properties of the projection operator $Q_{a_1,\cdots,a_{k-1}},$
\begin{eqnarray}\label{last}
f_{a_1,\cdots,a_k}(z)&=&f_{a_1,\cdots,a_{k-1}}(z)-\langle f_{a_1,\cdots,a_{k-1}},E_{a_1,\cdots,a_k}\rangle E_{a_1,\cdots,a_k}(z)\nonumber\\
&=&f_{a_1,\cdots,a_{k-1}}(z)-\left \langle f_{a_1,\cdots,a_{k-1}},B_{a_1,\cdots,a_{k-1}}e_{a_k}\right \rangle B_{a_1,\cdots,a_{k-1}}(z)e_{a_k}(z)\nonumber\\
&=&f_{a_1,\cdots,a_{k-1}}(z)-\left \langle \frac{f_{a_1,\cdots,a_{k-1}}}{B_{a_1,\cdots,a_{k-1}}},e_{a_k}
\right \rangle B_{a_1,\cdots,a_{k-1}}(z)e_{a_k}(z).
\end{eqnarray}
The modulus of (\ref{last}) is dominated by
\begin{eqnarray}\label{together}
|f_{a_1,\cdots,a_{k-1}}(z)|+\left |\frac{f_{a_1,\cdots,a_{k-1}}(a_k)}
{B_{a_1,\cdots,a_{k-1}}(a_k)}\right |
 \frac{1-|a_k|^2}{|1-\overline{a}_kz|}\end{eqnarray}
where the function
\[\frac{f_{a_1,\cdots,a_{k-1}}(z)}{B_{a_1,\cdots,a_{k-1}}(z)}\]
is analytic in $\overline{\mathbb{D}}.$ By invoking the Maximum Modulus Principle of analytic functions, it takes the maximal modulus on $\partial \mathbb{D}$ dominated by $3^{k-1}M,$ according to the inductive hypothesis. On the other hand, for $|z|\leq 1,$
\[ \left |\frac{1-|a_k|^2}{1-\overline{a}_kz}\right |\leq 2.\]
Altogether the quantity in (\ref{together}) is dominated by $3^kM,$ as desired.

\subsection{1-best approximation}
For $n=1, f\in \mathbb{H}^2,$ one can find $a_1\in \mathbb{D}$ such that  $|\langle f,E_{a_1}\rangle|=|f(a_1)|\sqrt{1-|a_1|^2}$ attains the maximal possible value of all the same kind. This is the so called Maximal Selection Principle proved through the boundary vanishing condition
\begin{eqnarray}\label{BVC} \lim_{|a|\to 1}|\langle f,E_{a}\rangle |=0\end{eqnarray}
via a Bolzano-Weierstrass compact argument\cite{QW1}. For the self-containing purpose we cite the simple proof of (\ref{BVC}) at the point.

For $\forall\epsilon>0$, we can find a polynomial function $g$ such that
\[ \|f-g\|<\frac{\epsilon}{2}.\]
Since $g$ is bounded in $\overline{\mathbb{D}},$ we have
\begin{eqnarray*} |\langle f,E_{a}\rangle|&\leq& |\langle g,E_{a}\rangle|+\epsilon/2\\
&=& \sqrt{1-|a|^2}|g(a)|+\epsilon/2\\
&\leq& \epsilon,
\end{eqnarray*}
if $|a|$ is sufficiently close to $1.$

 \subsection{2-best approximation}
By using the same density argument as for the $n=1$ case, we may assume that $f$ is a complex polynomial which is bounded, say by $M,$ in a neighbourhood of $\overline{\mathbb{D}}.$  Based on the definition of supreme, one can find a sequence of $2$-tuples, $(e_{a_1^{(l)}},e_{a_2^{(l)}}), l=1,2,\cdots,$ such that the norms of the projections $P_{{a_1^{(l)}},{a_2^{(l)}}}f$ tends to the supreme
 (\ref{sup}). Owing to continuity of inner product we may assume $a_1^{(l)}\ne a_2^{(l)}$ for every $l=1,2,\cdots.$ Since $(a_1^{(l)},a_2^{(l)})\in{\mathbb{D}}\times{\mathbb{D}},$ we may assume, through a Bolzano-Weierstrass compact argument, that the $2$-tuples $(a_1^{(l)},a_2^{(l)})$ converge to $(a_1,a_2)\in \overline{\mathbb{D}}\times \overline{\mathbb{D}}.$ If we can show $(a_1,a_2)\in {\mathbb{D}}\times {\mathbb{D}},$ then we are done. We show this by contradiction. Assume the opposite, that is, at least one of $a_1$ and $a_2$ is on the boundary $\partial \mathbb{D}.$ Since the projections $P_{{a_1^{(l)}},{a_2^{(l)}}}f$ are irrelevant with the order of ${a_1^{(l)}}, {a_2^{(l)}}$, we may assume that $a_2\in \partial \mathbb{D},$ and will then derive a contradiction.

  As a consequence of Lemma 1, for any $a_1$ in $\mathbb{D},$ there holds
$|f_{a_1}(z)|\leq 3M, z\in \mathbb{D}.$ By setting $f_1^{(l)}=f_{a^{(l)}_1}/(3M),$ we have $f_1^{(l)}(\mathbb{D})\subset \overline{\mathbb{D}}.$ Since $a^{(l)}_1\ne a^{(l)}_2$ and $f_{a^{(l)}_1}(a^{(l)}_1)=0$, a similar reasoning as in Lemma 1, we have
 \begin{eqnarray*}\label{coefficient2}
  |\langle f,E_{a^{(l)}_1,a^{(l)}_2}\rangle|
  &=& |\langle f_{a^{(l)}_1},E_{a^{(l)}_1,a^{(l)}_2}\rangle|\\
  &=& \left |\left \langle\frac{f_{a^{(l)}_1}}{B_{a^{(l)}_1}},e_{a^{(l)}_2}\right \rangle\right |\nonumber\\
  &=&\left |\frac{f_{a^{(l)}_1}(a^{(l)}_2)}{B_{a^{(l)}_1}(a^{(l)}_2)}\right |
  \sqrt{1-|a^{(l)}_2|^2}\nonumber\\
  &=&3M\left |\frac{f_1^{(l)}(a^{(l)}_2)}{B_{a^{(l)}_1}(a^{(l)}_2)}\right |
  \sqrt{1-|a^{(l)}_2|^2}\nonumber\\
  &=&3M\left |\frac{\frac{f_1^{(l)}(a^{(l)}_2)-f_1^{(l)}
  (a^{(l)}_1)}{1-\overline {f_1^{(l)}(a^{(l)}_1)}
  f_1^{(l)}(a^{(l)}_2)}}{\frac{a^{(l)}_2-a^{(l)}_1}
  {1-\overline{a^{(l)}_{1}}a^{(l)}_{2}}}\right |\sqrt{1-|a^{(l)}_2|^2}\nonumber\\
  &=&3M\frac{\rho(f_1^{(l)}(a^{(l)}_1),f_1^{(l)}(a^{(l)}_2))}
  {\rho(a^{(l)}_1,a^{(l)}_2)}\sqrt{1-|a^{(l)}_2|^2}.
\end{eqnarray*}
 Hence, the Schwarz lemma may be used to assert the boundedness of the first factor of the last term of the above chain of inequalities.  When $|a^{(l)}_2|\to 1,$ for $a^{(l)}_1$ uniformly,
 \begin{eqnarray*}\label{coefficient3}
  |\langle f,E_{a^{(l)}_1,a^{(l)}_2}\rangle|
  &=&3M\frac{\rho(f_1^{(l)}(a^{(l)}_1),f_1^{(l)}(a^{(l)}_2))}{\rho(a^{(l)}_1,a^{(l)}_2)}
  \sqrt{1-|a^{(l)}_2|^2}\\
  &\leq& 3M\sqrt{1-|a^{(l)}_2|^2}\to 0.
\end{eqnarray*}
Referring to (\ref{sup}), the above argument shows that $a_2^{(l)}$ does not help to get any larger $P_{a^{(l)}_1,a^{(l)}_2}f$ than $P_{a^{(l)}_1}f$, as $a_2^{(l)}$ tends to the boundary $\partial\mathbb{D}$.  This happens only in the case when $f$ is an $m$-Blaschke form with $m<2.$ In our case m=1 if $f$ is non-trivial, contradictory with the assumption of the theorem.

\subsection{General n-best approximation}
For a general $n$, as for the $n=2$ case, we are assuming that
 $f$ is a polynomial bounded by $M$ in a neighbourhood of $\overline{\mathbb{D}}.$ An analogous argument leads to a sequence of $n$-tuples $(a_1^{(l)},\cdots,a_{n}^{(l)})$ with mutually different terms that leads
 to
\begin{eqnarray}\label{become} \lim_{l\to \infty}\|P_{a_1^{(l)},\cdots,a_{n}^{(l)}}f\|=\sup \{\sum_{k=1}^n|\langle f,E_{b_1,b_2,...,b_k}\rangle|^2 \ :\ b_1,\cdots,b_n \in \mathbb{D}\}. \end{eqnarray}
Through a compact argument we may assume that the $n$-tuples $(a_1^{(l)},\cdots,a_{n}^{(l)})$ itself has a limit as an $n$-tuple $(a_1,\cdots,a_{n}),$ where the $a_k$'s are not necessarily mutually different and can be inside $\mathbb{D}$ or on the boundary of $\mathbb{D}$. If all the $a_k$'s are inside of $\mathbb{D},$ then
(\ref{become}) becomes
\begin{eqnarray}\label{be} \|P_{a_1,\cdots,a_{n}}f\|=\sup \{\sum_{k=1}^n|\langle f,E_{b_1,b_2,...,b_k}\rangle|^2 \ :\ b_1,\cdots,b_n \in \mathbb{D}\}, \end{eqnarray}
and we thus have the existence. We now show that this is indeed the case. We prove it by assuming the opposite and then derive a contradiction. Assume that at least one of the components sequences, say, $a^{(l)}_{k'}, 1\leq k'\leq n, l=1,2,\cdots,$ tends to the boundary $\partial\mathbb{D}.$ Since for each fixed $l,$ the projections $P_{a^{(l)}_1,\cdots,a^{(l)}_{n}}f$ are irrelevant with the order of $a^{(l)}_1,\cdots,a^{(l)}_{n}$, we may assume that $k'=n.$
Likewise to the $n=2$ case, set $$f^{(l)}_{n-1}(z)=\frac{f_{a^{(l)}_1,\cdots,a^{(l)}_{n-1}}(z)}{3^{n-1}M B_{a^{(l)}_1,a^{(l)}_2,...,a^{(l)}_{n-2}}(z)}.$$ Through analysis on the zeros of the denominator and the numerator functions, and invoking the Maximum Modulus Principle, this function is analytic for $\mathbb{D}\to\overline{\mathbb{D}}.$
 We have, when $|a^{(l)}_n|\to 1,$ for $a^{(l)}_1,a^{(l)}_2,...,a^{(l)}_{n-1}$ uniformly,
\begin{eqnarray*}
  |\langle f,E_{a^{(l)}_1,a^{(l)}_2,...,a^{(l)}_n}\rangle| &=& \left |\left \langle\frac{f_{a^{(l)}_1,a^{(l)}_2,...,a^{(l)}_{n-1}}}
  {B_{a^{(l)}_1,a^{(l)}_2,...,a^{(l)}_{n-1}}},e_{a^{(l)}_n}\right \rangle\right | \\
  &=&\left |\frac{f_{a^{(l)}_1,a^{(l)}_2,...,a^{(l)}_{n-1}}(a^{(l)}_n)}
  {B_{a^{(l)}_1,a_2,...,a^{(l)}_{n-1}}(a^{(l)}_n)}\right |\sqrt{1-|a^{(l)}_n|^2}\\
  &=&3^{n-1}M\left |\frac{\frac{f^{(l)}_{n-1}(a^{(l)}_n)-f^{(l)}_{n-1}(a^{(l)}_{n-1})}
  {1-\overline{f^{(l)}_{n-1}(a^{(l)}_n)}f^{(l)}_{n-1}(a^{(l)}_{n-1})}}
  {\frac{a^{(l)}_n-a^{(l)}_{n-1}}{1-\overline{a^{(l)}_{n-1}}a^{(l)}_n}}\right |\sqrt{1-|a^{(l)}_n|^2}\\
  &=&3^{n-1}M\frac{\rho(f^{(l)}_{n-1}(a^{(l)}_{n-1}),f^{(l)}_{n-1}(a^{(l)}_{n}))}
  {\rho(a^{(l)}_{n-1},a^{(l)}_n)}\sqrt{1-|a^{(l)}_n|^2}\\
  &\to& 0.
\end{eqnarray*}

\end{document}